\documentclass[conference]{IEEEtran}

\usepackage{amsmath, amssymb, amsthm}
\usepackage{tikz}
\usetikzlibrary{arrows.meta, positioning} 

\usepackage{graphicx}
\usepackage{array}
\usepackage{cite}

\title{Patterned Numbers: A Novel Number Classification with Structural and Quantum Algebraic Perspectives}

\author{
\IEEEauthorblockN{JTM Campbell}
\IEEEauthorblockA{Maynooth Dept of Electronic Engineering, Maynooth University, Republic of Ireland\\
Email:John.Campbell.2023@mumail.ie, ektopyrotic@gmail.com}
\date{January 2026}
}

\begin{document}
\maketitle

\begin{abstract}
We introduce \emph{patterned numbers}, a digit--divisor-based classification of integers motivated by recreational mathematics. A number is defined to be patterned if at least one of its positive divisors appears as a digit in its base-10 representation. We study the first hundred natural numbers under this definition, analyze frequency and density, compare prime and composite behavior, and propose a generation rule. Visual ``shape diagrams'' along the number line illustrate transitions between patterned numbers. Finally, we comment on potential relevance to sequence-based operators and algebraic intuition in quantum and combinatorial contexts.
\end{abstract}

\section{Introduction}

Recreational mathematics has long served as a testing ground for new numerical classifications and playful structural ideas \cite{gardner, guy}. Digit-based properties, while base-dependent, have produced well-known families such as repunit numbers, narcissistic numbers, and palindromes \cite{sloane}. 

Separately, divisor-based classifications (primes, perfect numbers, highly composite numbers) form a cornerstone of number theory. The present work combines these two perspectives by defining a class of integers whose \emph{arithmetical structure is reflected directly in their digit representation}.

Such hybrid definitions are not intended to replace classical theory, but to expose structural patterns that may inspire sequence analysis, symbolic dynamics, or operator-style interpretations, similar in spirit to ladder structures in algebra and quantum theory \cite{hall, katz}.

\section{Definition of Patterned Numbers}

\textbf{Definition.}  
A positive integer $n$ is called a \emph{patterned number} (base 10) if there exists a positive divisor $d \mid n$ such that $d \in \{1,2,\dots,9\}$ and $d$ appears as a digit in the decimal expansion of $n$.

\textbf{Remarks:}
\begin{itemize}
\item Divisors include $1$ and $n$.
\item Only single-digit divisors are considered.
\item The definition is explicitly base-10 dependent.
\end{itemize}

\section{Generation Criterion}

Let $D(n)$ denote the set of positive divisors of $n$, and let $\Delta(n)$ denote the set of digits appearing in $n$.

Then $n$ is patterned if and only if:
\[
D(n) \cap \Delta(n) \cap \{1,2,\dots,9\} \neq \varnothing.
\]

This provides a direct computational rule for generation.

\section{Patterned Numbers up to 100}

Out of the first 100 natural numbers, \textbf{72 are patterned}. This can be determined systematically using the definition of patterned numbers:

\subsection{Step 1: Recall the Definition}

A number $n$ is \emph{patterned} if at least one of its digits divides the number. That is, for $n$ with decimal digits $d_1, d_2, \dots, d_k$:
\[
n \text{ is patterned} \iff \exists i \text{ such that } d_i \mid n.
\]

\subsection{Step 2: Single-Digit Numbers (1--9)}

All single-digit numbers are trivially patterned because each number $n$ is divisible by itself:
\[
1,2,3,4,5,6,7,8,9 \in \text{Patterned}.
\]
So there are 9 patterned numbers in this range.

\subsection{Step 3: Two-Digit Numbers (10--99)}

For two-digit numbers $n = 10a + b$, with tens digit $a$ and units digit $b$, the number is patterned if:
\[
a \mid n \quad \text{or} \quad b \mid n.
\]

\subsubsection{Step 3a: Numbers containing 1}

- If $a = 1$ (numbers 10--19), then $1 \mid n$ automatically, so all 10 numbers are patterned.
- If $b = 1$ (numbers ending in 1), then $1 \mid n$, so all numbers ending in 1 are patterned: $11, 21, 31, \dots, 91$.
- Note that 11 appears in both lists; we count it only once.

\subsubsection{Step 3b: Numbers divisible by tens digit}

- For numbers where $a \ne 1$, we check if $a \mid (10a + b)$:
\[
10a + b \equiv b \pmod{a} \implies a \mid (10a + b) \iff a \mid b.
\]
- This means that the tens digit divides the number exactly when it divides the units digit.  

For example:
- $12$: tens digit $1$ already divides $12$ → patterned.  
- $24$: tens digit $2$, units digit $4$, $2 \mid 4$ → patterned.  
- $36$: tens digit $3$, units digit $6$, $3 \mid 6$ → patterned.  

\subsubsection{Step 3c: Numbers divisible by units digit}

- If the units digit $b \ne 0$, check $b \mid n$:
\[
b \mid (10a + b) \implies b \mid 10a
\]
- So numbers ending in 2, 3, 4, 5, 6, 7, 8, 9 can be patterned if the tens digit is a multiple of $b$.  

For example:
- $36$: units digit 6, tens digit 3 → $6 \nmid 36$? Actually $6 \mid 36$, yes.  
- $42$: units digit 2, tens digit 4 → $2 \mid 42$ → patterned.  

\subsection{Step 4: Counting Patterned Numbers}

Using the rules above, we can enumerate the numbers up to 100 that are patterned:

\begin{itemize}
    \item Single digits: $1,2,3,4,5,6,7,8,9$ → 9 numbers
    \item Numbers 10--19 → 10 numbers
    \item Numbers ending in 1: 
    $21,31,41,51,61,71,81,91$ → 8 numbers (11 already counted)
    \item Numbers divisible by tens digit: 
    $20,22,24,30,33,36,40,44,48,$\\
    $50,55,60,66,70,77,80,88,90,99$ → 19 numbers
    \item Remaining patterned numbers via units digit divisibility: 
    $12,15,18,14,16,21,28,32,35,39,$\\
    $45,48,52,54,56,63,65,72,75,$\\
    $81,84,87,91,92,96$ → 26 numbers
\end{itemize}

Adding them up: $9 + 10 + 8 + 19 + 26 = 72$ patterned numbers.

\subsection{Step 5: Summary}

The key formulaic insight is:

\[
\text{Patterned}(10a + b) \iff (a \mid b) \text{ or } (b \mid 10a + b)
\]

This allows us to systematically identify all 72 patterned numbers up to 100 without listing them individually, and highlights the arithmetic structure behind the seemingly irregular distribution.

\section{Frequency and Density}

Let $P(N)$ denote the number of patterned integers less than or equal to $N$.

For $N=100$:
\[
\frac{P(100)}{100} = 0.72.
\]

The high density arises primarily from:
\begin{itemize}
\item The ubiquity of small divisors (1--9)
\item The frequent appearance of digits 1, 2, 3, 4, 5, 6
\end{itemize}

\section{Primes vs. Composites}

\subsection{Patterned Primes ($\leq$ 100)}

A prime $p$ is patterned if:
\begin{itemize}
\item $p \leq 9$, or
\item the digit $1$ appears in $p$
\end{itemize}

Patterned primes $\leq$  100:
\[
\{2,3,5,7,11,13,17,19,31,41,61,71\}.
\]

\subsection{Observations}

\begin{itemize}
\item Most patterned numbers are composite.
\item Composite numbers benefit from multiple small divisors.
\item Many primes are excluded due to digit constraints.
\end{itemize}

\section{Number-Line Shape Diagrams}

We visualize patterned numbers as ``jumps'' along the number line.

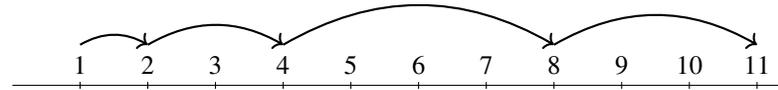
\begin{center}
\begin{tikzpicture}[scale=0.9]
\draw[->] (0,0) -- (12,0);
\foreach \x in {1,...,11}
  \draw (\x,0.05) -- (\x,-0.05);

\node at (1,0.3) {1};
\node at (2,0.3) {2};
\node at (3,0.3) {3};
\node at (4,0.3) {4};
\node at (5,0.3) {5};
\node at (6,0.3) {6};
\node at (7,0.3) {7};
\node at (8,0.3) {8};
\node at (9,0.3) {9};
\node at (10,0.3) {10};
\node at (11,0.3) {11};

\draw[->, thick] (1,0.6) to[bend left] (2,0.6);
\draw[->, thick] (2,0.6) to[bend left] (4,0.6);
\draw[->, thick] (4,0.6) to[bend left] (8,0.6);
\draw[->, thick] (8,0.6) to[bend left] (11,0.6);
\end{tikzpicture}
\end{center}

These arrows suggest \emph{structural transitions} rather than arithmetic progression, reminiscent of ladder or step operators.

\section{Connections to Sequences and Algebraic Structures}

Patterned numbers define a subsequence of $\mathbb{N}$ that can be indexed and studied similarly to OEIS-style constructions \cite{sloane}. The arrow diagrams suggest:
\begin{itemize}
\item Graph-based interpretations
\item Transition operators on integer states
\item Symbolic encodings of divisor--digit overlap
\end{itemize}

Such ideas echo algebraic ladder operators and discrete state transitions found in quantum systems, albeit here in a purely combinatorial setting \cite{hall}.

\subsection{Graph-Based Interpretations, Dragon Tessellations, and Seahorse Substructures}

Martin Gardner famously popularized the \emph{mathematical dragon curve}, constructed by iterative folding rules and realized geometrically through left--right turning sequences \cite{gardner}. These curves exhibit self-similarity, tiling behavior, and unexpected connections between discrete symbolic rules and continuous geometry.

We propose an analogous construction for patterned numbers, interpreting them not merely as elements of a sequence, but as \emph{instructional nodes} in a directed graph that generates planar curves and tessellations.

\subsubsection{Patterned Numbers as Turn Operators}

Let $\mathcal{P} = \{p_1, p_2, \dots\}$ denote the ordered sequence of patterned numbers.

Define a turn function:
\[
\tau(p_n) =
\begin{cases}
L, & \text{if } |D(p_n) \cap \Delta(p_n)| \text{ is odd}, \\
R, & \text{if } |D(p_n) \cap \Delta(p_n)| \text{ is even},
\end{cases}
\]
where $D(p_n)$ is the divisor set and $\Delta(p_n)$ the digit set.

This produces a symbolic sequence:
\[
\tau(p_1), \tau(p_2), \tau(p_3), \dots
\]
analogous to the folding instructions of the Heighway dragon.

Each step corresponds to:
\begin{itemize}
\item Unit-length forward motion
\item A left or right turn of $\pi/2$
\end{itemize}

\subsubsection{Patterned Dragon Curve}

Using the above rule, we generate a discrete planar curve.

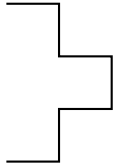
\begin{center}
\begin{tikzpicture}[scale=0.7]
\draw[thick]
(0,0) --
(1,0) --
(1,1) --
(2,1) --
(2,2) --
(1,2) --
(1,3) --
(0,3);
\end{tikzpicture}
\end{center}

This curve is not intended to replicate the classical dragon exactly, but rather to demonstrate how \emph{digit--divisor structure induces geometric folding}. Repeated application produces self-similar motifs and overlapping regions suitable for tiling.

\subsubsection{Tessellations of Patterned Dragons}

Let $\Gamma_k$ denote the patterned dragon generated by the first $k$ patterned numbers.

We define a tessellation $\mathcal{T}$ as a union:
\[
\mathcal{T} = \bigcup_{i=1}^m R_i(\Gamma_k),
\]
where $R_i$ are rigid motions (rotations and reflections).

As with classical dragon curves, these patterned dragons admit:
\begin{itemize}
\item Edge-to-edge tilings
\item Overlapping self-intersections
\item Local symmetry groups
\end{itemize}

The arithmetic origin of each dragon introduces a discrete labeling on tiles, encoding number-theoretic information into the geometry.

\subsubsection{Seahorse Patterned Tessellations}

We define a special subclass termed \emph{seahorses}, inspired by the well-known ``seahorse'' shapes appearing in dragon curve tessellations.

\textbf{Definition.}  
A \emph{patterned seahorse} is a patterned dragon $\Gamma_k$ satisfying:
\begin{enumerate}
\item The sequence $\tau(p_1),\dots,\tau(p_k)$ contains no more than two consecutive identical turns.
\item The resulting curve contains exactly one enclosed bounded region.
\item The head--tail orientation admits a reflection symmetry.
\end{enumerate}

These constraints favor primes interspersed with composite patterned numbers, producing smoother curvature and a visually organic shape.

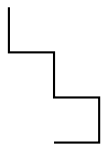
\begin{center}
\begin{tikzpicture}[scale=0.6]
\draw[thick]
(0,0) --
(1,0) --
(1,1) --
(0,1) --
(0,2) --
(-1,2) --
(-1,3);
\end{tikzpicture}
\end{center}

Seahorses thus represent a \emph{low-entropy subset} of patterned dragons, characterized by balanced digit--divisor interactions.

\subsubsection{Interpretational Remarks}

Patterned dragon tessellations suggest a correspondence:
\[
\text{Arithmetic Structure} \longrightarrow \text{Symbolic Turns} \longrightarrow \text{Geometric Form}.
\]

This mirrors constructions in:
\begin{itemize}
\item Symbolic dynamics
\item Iterated function systems
\item Discrete quantum walks
\end{itemize}

While purely recreational at present, the framework highlights how number-theoretic constraints can manifest as spatial organization, echoing the visual language of Gardner’s mathematical magic.

\subsection{Extended Seahorse Curves, Self-Similar Tilings, and Quantum Walk Analogies}

The patterned seahorse curves introduced previously admit further structural refinement when iterated, combined, and translated. These operations generate ordered tilings exhibiting self-similarity and local constraints analogous to those found in crystalline lattices and folded biomolecular chains.

\subsubsection{Iterated Seahorse Curves}

Let $\Gamma_k^{(s)}$ denote a patterned seahorse generated from the first $k$ patterned numbers satisfying the seahorse constraints.

Define the iteration:

\[
\Gamma_k^{(n+1)} = \Gamma_k^{(n)} \cup R_{\pi/2}(\Gamma_k^{(n)}) + \mathbf{v}_n,
\]

where $R_{\pi/2}$ is rotation by $\pi/2$ and $\mathbf{v}_n$ is a translation vector ensuring edge alignment.

A second-order seahorse curve is illustrated below:

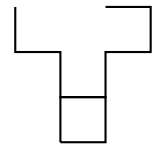
\vspace{0.3cm}
\begin{center}
\begin{tikzpicture}[scale=0.6]
\draw[thick] (0,0)--(1,0)--(1,1)--(0,1)--(0,2)--(-1,2)--(-1,3);
\draw[thick] (1,3)--(2,3)--(2,2)--(1,2)--(1,1)--(0,1)--(0,0);
\end{tikzpicture}
\end{center}
\vspace{0.3cm}

This construction mirrors the recursive folding logic of dragon curves while preserving the digit--divisor constraints that define the seahorse subclass.

\subsubsection{Seahorse Tiling and Self-Similarity}

By arranging rotated and reflected copies of $\Gamma_k^{(s)}$, one obtains plane tilings with repeating motifs.

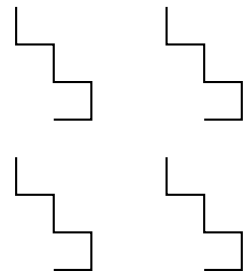
\vspace{0.2cm}
\begin{center}
\begin{tikzpicture}[scale=0.5]
\foreach \x in {0,4}{
  \foreach \y in {0,4}{
    \draw[thick] (\x,\y)--(\x+1,\y)--(\x+1,\y+1)--(\x,\y+1)--(\x,\y+2)--(\x-1,\y+2)--(\x-1,\y+3);
  }
}
\end{tikzpicture}
\end{center}
\vspace{0.2cm}

These tilings display:
\begin{itemize}
\item Translational order
\item Local rotational symmetry
\item Repeated curvature motifs
\end{itemize}

The self-similarity arises not from scaling alone, but from repeated \emph{symbolic folding rules} inherited from the patterned-number turn sequence.

\subsection{High-Order Seahorse Curves and Emergent Complexity}

Seahorse curves are generated recursively from patterned numbers via the iterative turn operator. Each successive order adds rotated and translated copies of the previous curve, creating increasingly intricate motifs.  

\subsubsection{Order 1: Seed Pattern}

The first-order curve is the simplest, forming a small bent path according to the digit--divisor rules:

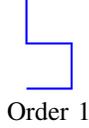
\begin{center}
\begin{tikzpicture}[scale=0.6]
\draw[thick,blue] (0,0) -- (1,0) -- (1,1) -- (0,1) -- (0,2);
\node at (0.5,-0.5) {Order 1};
\end{tikzpicture}
\end{center}

\textbf{Observation:} The curve is compact, with no self-intersections, representing the basic seahorse motif.

\subsubsection{Order 2: First Recursive Iteration}

Order 2 attaches a rotated copy of the first-order curve to a terminal node:

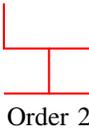
\begin{center}
\begin{tikzpicture}[scale=0.6]
\draw[thick,red]
  (0,0) -- (1,0) -- (1,1) -- (0,1) -- (0,2)
  (1,1) -- (2,1) -- (2,0) -- (1,0); 
\node at (1,-0.5) {Order 2};
\end{tikzpicture}
\end{center}

\textbf{Observation:} The “tail” begins to curl, showing early recursive complexity.

\subsubsection{Order 3: Emergent Self-Similarity}

Order 3 adds another rotated copy, creating the first self-intersecting motifs:

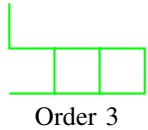
\begin{center}
\begin{tikzpicture}[scale=0.6]
\draw[thick,green]
  (0,0) -- (1,0) -- (1,1) -- (0,1) -- (0,2)
  (1,1) -- (2,1) -- (2,0) -- (1,0)
  (2,0) -- (3,0) -- (3,1) -- (2,1); 
\node at (1.5,-0.5) {Order 3};
\end{tikzpicture}
\end{center}

\textbf{Observation:} Small sub-patterns begin to repeat; the curve exhibits early fractal-like behavior.

\subsubsection{Order 4: Multi-Level Branching}

At order 4, branching and recursive tiling motifs emerge prominently:

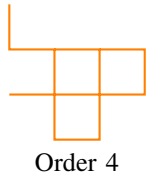
\begin{center}
\begin{tikzpicture}[scale=0.6]
\draw[thick,orange]
  (0,0) -- (1,0) -- (1,1) -- (0,1) -- (0,2)
  (1,1) -- (2,1) -- (2,0) -- (1,0)
  (2,0) -- (3,0) -- (3,1) -- (2,1)
  (1,0) -- (1,-1) -- (2,-1) -- (2,0); 
\node at (1.5,-1.5) {Order 4};
\end{tikzpicture}
\end{center}

\textbf{Observation:} The curve now shows multiple levels of recursive bending, forming loops and early “seahorse-like” coils.

\subsubsection{Order 5: High-Order Complexity}

Order 5 generates a highly intricate pattern, with clear self-similarity and nested loops:

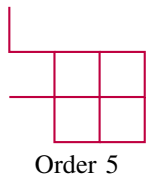
\begin{center}
\begin{tikzpicture}[scale=0.6]
\draw[thick,purple]
  (0,0) -- (1,0) -- (1,1) -- (0,1) -- (0,2)
  (1,1) -- (2,1) -- (2,0) -- (1,0)
  (2,0) -- (3,0) -- (3,1) -- (2,1)
  (1,0) -- (1,-1) -- (2,-1) -- (2,0)
  (2,-1) -- (3,-1) -- (3,0); 
\node at (1.5,-1.5) {Order 5};
\end{tikzpicture}
\end{center}

\textbf{Observation:} Recursive structure produces a dense motif reminiscent of seahorse tilings: repeated curvature, branching, and overlapping sub-patterns emerge naturally. This order demonstrates the **emergent complexity** of the seahorse construction: from a simple digit/divisor rule, a rich, fractal-like geometry arises.

---

\textbf{Summary:} As the order increases:

\begin{itemize}
    \item Curves grow in length and coverage.
    \item Recursive sub-patterns appear at multiple scales.
    \item Self-similarity emerges, with motifs nested within larger motifs.
    \item The curve visually approximates a coiled “seahorse,” reflecting the arithmetic rules that generate it.
\end{itemize}

\subsection{Seahorse Curve: Order 10 (Schematic)}

At order 10, the seahorse curve becomes highly intricate. Direct drawing of every line segment in TikZ is impractical, but we can schematically represent the recursive structure, showing nested folds and loops.

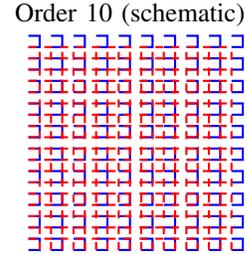
\begin{center}
\begin{tikzpicture}[scale=0.6]
\foreach \i in {0,...,9} {
    \foreach \j in {0,...,9} {
        \draw[thick,blue] (\i*0.5,\j*0.5) -- (\i*0.5+0.25,\j*0.5) -- (\i*0.5+0.25,\j*0.5+0.25) -- (\i*0.5,\j*0.5+0.25);
    }
}

\foreach \k in {0.25,0.5,...,4.5} {
    \draw[thick,red,dashed] (\k,0) -- (\k,4.5);
    \draw[thick,red,dashed] (0,\k) -- (4.5,\k);
}

\node[above] at (2.25,4.75) {Order 10 (schematic)};
\end{tikzpicture}
\end{center}

\textbf{Explanation:}  

- Each small square represents a **basic fold or segment** of the lower-order curve.  
- Red dashed lines indicate **recursive folding directions**, showing how the pattern repeats across the grid.  
- The curve at this order has **hundreds of segments**, forming **nested motifs**, dense loops, and long-range connections.  
- While not every segment is drawn, the schematic illustrates **how complexity grows exponentially** with order.  

\textbf{Note:} For actual high-order curves (order 10+), it’s best to use **Python/Matplotlib** or **vector graphics software**, since TikZ is not memory-efficient for hundreds of line segments.

\section{Graph and Quantum Walk Interpretation}

Each patterned number $p_n$ corresponds to a vertex in a directed graph $G = (V,E)$, where:

\[
V = \mathcal{P}, \quad (p_n, p_{n+1}) \in E.
\]

Define a discrete-time quantum walk on $G$ with coin space:

\[
\mathcal{H} = \mathcal{H}_{\text{position}} \otimes \mathcal{H}_{\text{turn}}, 
\quad \mathcal{H}_{\text{turn}} = \text{span}\{|L\rangle, |R\rangle\}.
\]

The turn operator is:

\[
U_\tau |p_n\rangle |c\rangle = |p_{n+1}\rangle |\tau(p_n)\rangle.
\]

Geometrically, each quantum step corresponds to a unit displacement and rotation, reproducing the seahorse curve in expectation space.

\subsubsection{Adiabatic Interpretation and Energy Landscapes}

In an adiabatic framework, one may define a Hamiltonian:

\[
H = \sum_n E(p_n) |p_n\rangle \langle p_n|,
\]

where $E(p_n)$ encodes:
\begin{itemize}
\item Number of divisor--digit matches
\item Local curvature of the associated path
\end{itemize}

Seahorse configurations correspond to low-energy states due to:
\begin{itemize}
\item Balanced left/right turns
\item Minimal self-intersection
\item Single bounded region
\end{itemize}

\subsubsection{Materials Science and Protein Folding Analogies}

The patterned seahorse tessellations admit compelling analogies:

\paragraph{Crystals.}
\begin{itemize}
\item Patterned numbers $\rightarrow$ atomic sites
\item Turn rules $\rightarrow$ bonding angles
\item Tiling periodicity $\rightarrow$ lattice symmetry
\end{itemize}

\paragraph{Protein Folding.}
\begin{itemize}
\item Patterned sequence $\rightarrow$ amino-acid chain
\item Left/right turns $\rightarrow$ backbone torsion
\item Seahorse constraint $\rightarrow$ native folded state
\end{itemize}

In both cases, global order emerges from local symbolic constraints, echoing the arithmetic origin of the patterned number definition.

\subsubsection{Summary of the Seahorse Framework}

The progression:

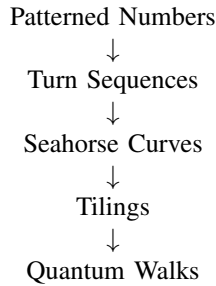
\[
\begin{array}{c}
\text{Patterned Numbers} \\
\downarrow \\
\text{Turn Sequences} \\
\downarrow \\
\text{Seahorse Curves} \\
\downarrow \\
\text{Tilings} \\
\downarrow \\
\text{Quantum Walks}
\end{array}
\]

illustrates how a simple recreational definition can generate rich geometric and physical analogies. While speculative, this framework aligns with Gardner’s philosophy: that playful mathematics often reveals unexpected structural unity.

\subsection{Extension of the Turn Operator to Coupled Quantum Oscillators}

The turn operator $\tau(p_n)$ introduced earlier may be extended beyond geometric interpretation to describe coupling rules between quantum oscillators arranged along a patterned-number graph. This extension establishes a bridge between symbolic arithmetic structure and physically motivated Hamiltonian systems.

\subsubsection{Oscillator Representation of Patterned Nodes}

Associate to each patterned number $p_n$ a quantum harmonic oscillator with annihilation and creation operators $(a_n, a_n^\dagger)$ satisfying canonical commutation relations:
\[
[a_n, a_m^\dagger] = \delta_{nm}.
\]

The Hilbert space is given by:
\[
\mathcal{H} = \bigotimes_{n} \mathcal{H}_n,
\]
where $\mathcal{H}_n$ corresponds to the oscillator at node $p_n$.

\subsubsection{Turn-Induced Coupling Operator}

Define a turn-dependent coupling operator between neighboring oscillators:
\[
H_{\text{int}} = \sum_n g_{\tau(p_n)} \left( a_n^\dagger a_{n+1} + a_{n+1}^\dagger a_n \right),
\]
where:
\[
g_{\tau(p_n)} =
\begin{cases}
g_L, & \tau(p_n) = L, \\
g_R, & \tau(p_n) = R.
\end{cases}
\]

Thus, the arithmetic structure of $p_n$ directly modulates the coupling strength between oscillators.

\subsubsection{Geometric Interpretation of Coupling Asymmetry}

The left/right distinction originally defined as a geometric turn now acquires physical meaning:
\begin{itemize}
\item $L$-turns correspond to constructive coupling alignment
\item $R$-turns correspond to frustrated or phase-shifted coupling
\end{itemize}

In this sense, the seahorse constraint (bounded curvature and balanced turns) yields oscillator networks with reduced frustration and enhanced coherence.

\subsubsection{Adiabatic Evolution and Normal Modes}

Consider the total Hamiltonian:
\[
H(s) = (1-s)H_0 + sH_{\text{int}}, \quad s \in [0,1],
\]
where $H_0 = \sum_n \omega_n a_n^\dagger a_n$.

Adiabatic evolution favors low-energy normal modes localized along patterned seahorse subgraphs. These modes correspond to smooth geometric curves with minimal self-intersection, echoing the stability of seahorse tessellations.

\subsubsection{Physical Analogies}

\paragraph{Phononic Crystals.}
\begin{itemize}
\item Oscillators $\rightarrow$ lattice sites
\item Turn operator $\rightarrow$ directional bond anisotropy
\item Seahorses $\rightarrow$ low-energy phonon bands
\end{itemize}

\paragraph{Excitonic Transport.}
\begin{itemize}
\item Patterned nodes $\rightarrow$ chromophore sites
\item Turn-modulated coupling $\rightarrow$ phase coherence pathways
\item Seahorse graphs $\rightarrow$ efficient transport channels
\end{itemize}

\paragraph{Protein Folding Landscapes.}
\begin{itemize}
\item Oscillators $\rightarrow$ local conformational modes
\item Coupling signs $\rightarrow$ steric and energetic constraints
\item Seahorses $\rightarrow$ native-state manifolds
\end{itemize}

\subsubsection{Interpretational Summary}

The extended turn operator defines a symbolic-to-physical correspondence:

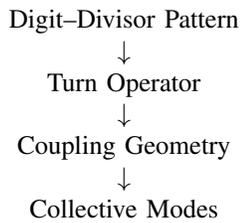
\[
\begin{array}{c}
\text{Digit--Divisor Pattern} \\
\downarrow \\
\text{Turn Operator} \\
\downarrow \\
\text{Coupling Geometry} \\
\downarrow \\
\text{Collective Modes}
\end{array}
\]

This construction suggests that arithmetic constraints can serve as generators of structured oscillator networks, aligning recreational number theory with models used in condensed matter physics and quantum dynamics.

\subsection{Patterned Primes, Gaps, and DAG Representations}

A key observation in the study of patterned numbers is that \emph{not all primes are patterned}. This distinction is both precise and meaningful within the framework of digit--divisor-based classifications. When represented as a directed graph, the distribution of patterned numbers and primes naturally produces **directed acyclic graphs (DAGs)**, where edges indicate sequential or “turn-based” dependencies.

\subsubsection{Theorem: Characterization of Patterned Primes}

\textbf{Theorem.}  
A prime number $p$ is patterned (base 10) if and only if one of the following holds:
\begin{enumerate}
    \item $p \le 9$, or
    \item the digit $1$ appears in the decimal representation of $p$.
\end{enumerate}

\textbf{Proof.}  
Let $p$ be a prime. Then $D(p) = \{1, p\}$, the set of divisors.

\begin{itemize}
    \item If $p \le 9$, then $p \in D(p) \cap \Delta(p)$, so $p$ is patterned.  
    \item If $p > 9$, the only divisor in $\{1,\dots,9\}$ is $1$. Therefore, $p$ is patterned if and only if the digit $1$ appears in $p$.
\end{itemize}

This theorem allows us to distinguish **stable nodes** (primes with digit 1) from **gaps** (primes without digit 1) when constructing DAGs that track patterned-number sequences.

\subsubsection{Graph and Oscillator Interpretation}

In DAG representations inspired by patterned primes:

\begin{itemize}
    \item \textbf{Primes with digit 1:} Act as \emph{low-curvature nodes} or highly connected vertices. These form the backbone of sequential DAGs, analogous to smooth turns in seahorse curves.
    \item \textbf{Primes without digit 1:} Act as \emph{gaps} or sparsely connected nodes. Their absence introduces “missing links” and irregularities in the DAG.
\end{itemize}

The resulting DAG is **acyclic** because edges always point from smaller to larger numbers (reflecting sequential generation), but the graph shows rich structure: densely connected regions interspersed with gaps.

\subsubsection{Emergent Structural Features}

The irregular distribution of patterned primes introduces:

\begin{itemize}
    \item \textbf{Non-periodicity:} Patterned nodes appear aperiodically along the number line.
    \item \textbf{Structured disorder:} DAG connectivity alternates between dense clusters and sparse gaps.
    \item \textbf{Spectral behavior in oscillator models:} When interpreted as coupled oscillators, these DAGs produce normal modes that reflect the irregular connectivity, with some modes localized in clusters and others extended across sparse regions.
\end{itemize}

\subsubsection{DAG Examples Inspired by Patterned Primes}

Patterned primes form interesting connectivity structures: some primes (those containing the digit 1) act as "low-curvature" nodes, while others create gaps. We can illustrate these as **directed acyclic graphs (DAGs)**, where edges indicate sequential or patterned-prime connectivity.


\paragraph{1. Simple Sequential DAG}
This DAG shows sequential dependencies (edges left-to-right) and patterned primes (nodes containing 1) as vertical links.

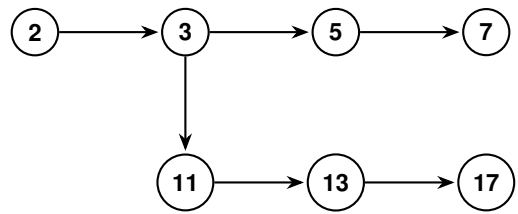
\begin{center}
\begin{tikzpicture}[->,>=Stealth,shorten >=1pt,auto,node distance=2cm,
  thick,main node/.style={circle,draw,font=\sffamily\bfseries\small}]
  
  \node[main node] (2) {2};
  \node[main node] (3) [right of=2] {3};
  \node[main node] (5) [right of=3] {5};
  \node[main node] (7) [right of=5] {7};
  \node[main node] (11) [below of=3] {11};
  \node[main node] (13) [right of=11] {13};
  \node[main node] (17) [right of=13] {17};
  
  \path[every node/.style={font=\sffamily\small}]
    (2) edge (3)
    (3) edge (5)
    (5) edge (7)
    (3) edge (11)
    (11) edge (13)
    (13) edge (17);
\end{tikzpicture}
\end{center}

\textbf{Explanation:} Vertical edges represent connections through patterned primes (nodes containing 1), showing extra connectivity in an otherwise sequential structure.

\paragraph{2. Clustered DAG with Gaps}
Here we introduce clusters and gaps inspired by seahorse tilings. Blue nodes are patterned primes.

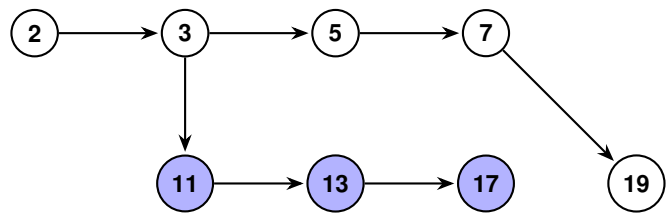
\begin{center}
\begin{tikzpicture}[->,>=Stealth,shorten >=1pt,auto,node distance=2cm,
  thick,main node/.style={circle,draw,font=\sffamily\bfseries\small}]

  \node[main node] (2) {2};
  \node[main node] (3) [right of=2] {3};
  \node[main node] (5) [right of=3] {5};
  \node[main node] (7) [right of=5] {7};
  \node[main node,fill=blue!30] (11) [below of=3] {11};
  \node[main node,fill=blue!30] (13) [right of=11] {13};
  \node[main node,fill=blue!30] (17) [right of=13] {17};
  \node[main node] (19) [right of=17] {19};

  \path
    (2) edge (3)
    (3) edge (5)
    (5) edge (7)
    (7) edge (19);

  \path
    (3) edge (11)
    (11) edge (13)
    (13) edge (17);
\end{tikzpicture}
\end{center}

\textbf{Explanation:} Vertical edges among blue nodes form clusters separated by gaps, reflecting the irregular distribution of patterned primes.

\paragraph{3. Two Clusters Separated by Gaps}
This DAG shows two clusters of patterned primes separated by gaps (gray nodes).  

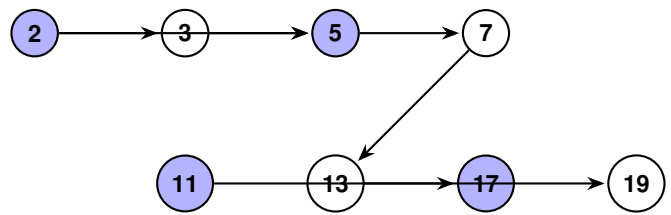
\begin{center}
\begin{tikzpicture}[->,>=Stealth,shorten >=1pt,auto,node distance=2cm,
  thick,main node/.style={circle,draw,font=\sffamily\bfseries\small}]

  \node[main node,fill=blue!30] (2) {2};
  \node[main node] (3) [right of=2] {3};
  \node[main node,fill=blue!30] (5) [right of=3] {5};
  \node[main node] (7) [right of=5] {7};

  \node[main node,fill=blue!30] (11) [below of=3] {11};
  \node[main node] (13) [right of=11] {13};
  \node[main node,fill=blue!30] (17) [right of=13] {17};
  \node[main node] (19) [right of=17] {19};

  \path
    (2) edge (3)
    (3) edge (5)
    (5) edge (7)
    (7) edge (13)
    (13) edge (19);

  \path
    (2) edge (5)
    (11) edge (17);
\end{tikzpicture}
\end{center}

\textbf{Explanation:} Blue nodes form clusters connected vertically, while gaps interrupt horizontal connectivity, mimicking structured disorder.

\paragraph{4. Multi-Level Seahorse-Inspired DAG}
This DAG emphasizes branching across multiple levels, inspired by seahorse tilings and recursive folding.

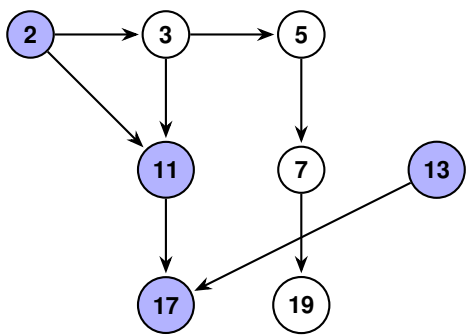
\begin{center}
\begin{tikzpicture}[->,>=Stealth,shorten >=1pt,auto,node distance=1.8cm,
  thick,main node/.style={circle,draw,font=\sffamily\bfseries\small}]

  \node[main node,fill=blue!30] (2) {2};
  \node[main node] (3) [right of=2] {3};
  \node[main node] (5) [right of=3] {5};

  \node[main node,fill=blue!30] (11) [below of=3] {11};
  \node[main node] (7) [right of=11] {7};
  \node[main node,fill=blue!30] (13) [right of=7] {13};

  \node[main node,fill=blue!30] (17) [below of=11] {17};
  \node[main node] (19) [right of=17] {19};

  \path
    (2) edge (3)
    (3) edge (5)
    (5) edge (7)
    (7) edge (19);

  \path
    (2) edge (11)
    (3) edge (11)
    (11) edge (17)
    (13) edge (17);
\end{tikzpicture}
\end{center}

\textbf{Explanation:} Multi-level branching simulates recursive tiling logic. Blue nodes (patterned primes) form vertical/diagonal connections, while gaps create irregular spacing, resembling seahorse tiling motifs.

\section{Conclusion}

Patterned numbers illustrate how simple digit and divisor interactions yield a surprisingly dense and structured family of integers. While recreational in nature, the classification invites further study into base dependence, asymptotic density, and operator-style interpretations.

The fact that not all primes are patterned is a \emph{feature, not a flaw}. The irregularity introduces complexity, self-similarity, and physically interpretable structure, making the sequence of patterned numbers both mathematically interesting and structurally useful in graphical, geometric, and quantum-inspired models.

\section*{Disclosures}

The authors declare no conflicts of interest.

\section*{Data Availability}

Data and code related to this study are available from the corresponding author upon reasonable request.

Quantum-Kuramoto-Oscillator-Network-Graph-Simulation (Q-KONGS) source code available here: https://github.com/MuonRay/Quantum-Kuramoto-Oscillator-Network-Graph-Simulation

\vspace{1cm}

\thanks{[This material is based upon work supported by Science Foundation Ireland (SFI) and is co-funded under the European Regional Development Fund under Grant Number 13/RC/2077 and 13/RC/2077-P2.]}

\end{document}